\documentclass[12pt]{amsart}
\usepackage{amssymb}
\usepackage[all]{xy}
\usepackage{fullpage}
\usepackage{diagbox}
\usepackage[pdfencoding=auto]{hyperref}
\allowdisplaybreaks

\numberwithin{equation}{section}

\theoremstyle{plain}
\newtheorem{theorem}[equation]{Theorem}
\newtheorem{lemma}[equation]{Lemma}
\newtheorem{proposition}[equation]{Proposition}
\newtheorem{cor}[equation]{Corollary}

\theoremstyle{definition}

\theoremstyle{remark} 
\newtheorem{remark}[equation]{Remark} 

\newtheorem{definition}[equation]{Definition}

\newcommand{\bF}{\mathbb{F}}

\newcommand{\bN}{\mathbb{N}}

\newcommand{\bZ}{\mathbb{Z}}

\newcommand{\cG}{\mathcal{G}}

\newcommand{\cP}{\mathcal{P}}

\renewcommand{\ge}{\geqslant}
\renewcommand{\geq}{\geqslant}
\newcommand{\HH}{H\!H}
\newcommand{\Hom}{\mathsf{Hom}}

\newcommand{\lcm}{\mathsf{lcm}}
\renewcommand{\le}{\leqslant}
\renewcommand{\leq}{\leqslant}

\newcommand{\sgn}{\mathsf{sgn}}
\newcommand{\sym}{\mathfrak S}

\author{David Benson, Radha Kessar, and Markus Linckelmann}

\address{David Benson \\ 
Institute of Mathematics\\ 
Fraser Noble Building\\
University of Aberdeen\\ 
King's College\\ 
Aberdeen AB24 3UE\\ 
United Kingdom}

\address{Radha Kessar and Markus Linckelmann \\
School of Mathematics, Computer Science \& Engineering \\
Department of Mathematics \\
City, University of London \\
Northampton Square \\
London EC1V 0HB \\
United Kingdom}



\title{Hochschild cohomology of symmetric groups in low degrees}

\begin{document}

\begin{abstract}
We compute the dimensions of the Hochschild cohomology of 
symmetric groups over prime fields in low degrees.
This involves us in studying some partition identities and generating
functions of the dimensions in any fixed degree of the Hochschild 
cohomology of symmetric groups.
\end{abstract}

\maketitle

\section{Introduction}

The purpose of this note is to compute the dimensions of the
Hochschild cohomology of symmetric groups over prime fields
in low degrees. We relate this to some partition identities.

For $n \ge 0$ we denote by $p(n)$ the number of partitions of $n$,
with the convention $p(0)=1$.
We write $P(t)$ for the generating function of $p(n)$, which is given by Euler's
identity as
\begin{equation}\label{eq:P(t)} 
P(t) = \sum_{n=0}^\infty p(n)\,t^n = \prod_{m=1}^\infty
  \frac{1}{1-t^m} = 1\bigg/ \sum_{n=-\infty}^\infty (-1)^n
  t^{n(3n+1)/2}. 
\end{equation}

Throughout the paper $p$ is a prime. Our main theorem  is as follows.

\begin{theorem}\label{th:HH}
The generating functions for the dimensions of the Hochschild cohomology of the group
algebra of the symmetric group $\sym_n$ on $n$ letters  in low degrees are given by
\begin{enumerate}
\item[\rm (i)]\qquad
$\displaystyle\sum_{n=0}^\infty \dim_{\bF_p}(\HH^0(\bF_p\sym_n))\,t^n=P(t)$
\item[\rm (ii)]\qquad
$\displaystyle\sum_{n=0}^\infty \dim_{\bF_p}(\HH^1(\bF_p\sym_n))\,t^n=
\begin{cases}
\displaystyle\frac{2t^2}{1-t^2}\,P(t) &p=2, \\[10pt]
\displaystyle\frac{t^p}{1-t^p}\,P(t) & p \ge 3.
\end{cases}$\medskip
\item[\rm (iii)]\qquad
$\displaystyle\sum_{n=0}^\infty \dim_{\bF_p}(\HH^2(\bF_p\sym_n))\,t^n=
\begin{cases}
\displaystyle\frac{2t^2+3t^4-t^6}{(1-t^2)(1-t^4)}\,P(t)
&p=2,\\[10pt]
\displaystyle\frac{t^p}{(1-t^p)(1-t^{2p})}\,P(t) & p \ge 3.
\end{cases} $
\end{enumerate}
\end{theorem}

For degree zero, this is well known: 
$p(n)$ is equal to the dimension of the centre of $\bF_p\sym_n$.
For degrees one and two,
the proof of the theorem is given in
Sections~\ref{se:HH1} and~\ref{se:HH2}, where we give explicit combinatorial
descriptions of the dimensions of $\HH^1(\bF_p \sym_n)$ and $\HH^2(\bF_p \sym_n)$
in Corollary \ref{cor:HH1} and Corollary \ref{cor:HH2}, respectively.

In general,  we have the  following  result.

\begin{theorem} \label{th:rational}
For any integer $r\geq 0$ there exists a rational function $R_{p,r}(t)$ 
such that 
$$\sum_{n = 0}^\infty \dim_{\bF_p}(\HH^r(\bF_p\sym_n))\,  t^n =
R_{p,r}(t) P(t).$$
\end{theorem}

The theorem is proved in Section~\ref{se:HHgeq3}. The  analysis  in    Section~\ref{se:HHgeq3}      could be used  to  describe   $R_{p,1}(t)$ and  $R_{p,2}(t)$ as in Theorem~\ref{th:HH} but we 
have chosen to treat these cases  separately for  the  sake of  
obtaining    explicit combinatorial descriptions of the dimensions of 
$\HH^1(\bF_p \sym_n)$ and $\HH^2(\bF_p \sym_n)$.

\bigskip\noindent
{\bf Acknowledgements.} 
The first author is grateful to City, University of London for its
hospitality during the research for this paper. The second author  acknowledges support from
EPSRC grant EP/T004592/1.

\section{Some partition identities}\label{se:partitions}  

In this section, we discuss some partition identities that will appear
in the proofs in later sections. For $n \ge 0$
we write $\cP(n)$ for the set of partitions of $n$, with the standard convention that
$\cP(0)$ contains the empty partition as unique element. We set
$p(n) =$  $|\cP(n)|$, the number of partitions of $n$, and we write
$P(t)$ for the generating function, as described in~\eqref{eq:P(t)}.
If $\lambda$ is a partition of $n$, we write $\lambda \vdash n$. 
We write $\lambda_k$ for the number of parts of $\lambda$ of length $k$,
and we use the notation $\lambda = (n^{\lambda_n} \dots 2^{\lambda_2} 1^{\lambda_1})$.

\begin{definition}\label{def:f}
Let $n$, $k$ be integers such that $n\ge 0$ and $k\ge 1$.
Define  $\displaystyle F_k(n)=\sum_{\lambda\vdash n}\lambda_k$. 
Thus $F_k(n)$ is
the total number of parts of length $k$ in all partitions of $n$.
\end{definition}

\begin{lemma}\label{le:F}
Let $n$, $k$ be  integers such that
$n\ge k\ge 1$. We have $F_k(n)=F_k(n-k)+p(n-k)$.
\end{lemma}
\begin{proof}
Removing a part of length $k$ gives a bijection between the partitions of $n$
with a part of length $k$ and the partitions of $n-k$, which reduces by one
the number of parts of length $k$.
\end{proof}

\begin{proposition}\label{pr:F}
Let $n$, $k$ be integers such that $n\ge 0$ and $k\ge 1$. We have
\[ F_k(n) = \sum_{i=1}^{\lfloor\frac{n}{k}\rfloor}p(n-ik)=
p(n-k) + p(n-2k)+ \cdots \]
In particular, $F_k(n)=0$ for $n<k$.
The generating function for $F_k(n)$ is given by
\[ \sum_{n=1}^\infty F_k(n)\,t^n = \frac{t^k}{1-t^k}\,P(t). \]
\end{proposition}
\begin{proof}
The first part follows by induction from Lemma~\ref{le:F}.
Thus 
\begin{align*} 
\sum_{n=1}^\infty F_k(n)\,t^n &= 
\sum_{n=1}^\infty\left(\sum_{i=1}^{\lfloor\frac{n}{k}\rfloor}
p(n-ik)\right)\,t^n 
=\sum_{i=1}^{\infty}\sum_{n=ik}^\infty p(n-ik)\,t^n \\
&= 
\sum_{i=1}^\infty\sum_{n=0}^\infty p(n)\,t^{n+ik} =
\sum_{i=1}^{\infty} P(t)\,t^{ik} = \frac{t^k}{1-t^k}\,P(t).
\qedhere
\end{align*}
\end{proof}

For $1\le k\le 10$, the sequences $F_k(n)$ can be found in Sloane's
online encyclopedia of integer sequences as sequence number 
A000070 for $k=1$, and
A024784${}+k$ for $2\le k\le 10$.\bigskip

Next, we relate $F_k(n)$ to some other functions defined in terms of partitions.
These are the functions that we associate with Hochschild cohomology
of symmetric groups, through the centraliser decomposition.

\begin{definition}
Let $n$ be a non-negative integer.
If $\lambda\vdash n$ and $k$, $\ell$ are positive integers,  let 
$g_{k,\ell}(\lambda)$ be the number of distinct lengths
divisible by $k$ of parts of $\lambda$ that are repeated at least $\ell$ times.
More precisely, if $\lambda=(n^{\lambda_n} \dots 2^{\lambda_2}1^{\lambda_1})$,
then $g_{k,\ell}(\lambda)$ is the number of indices $i$, $1\le i\le n$, 
such that $k$ divides $i$ and $\lambda_i\ge \ell$.
 Let 
$G_{k,\ell}(n)=\displaystyle\sum_{\lambda\vdash n} g_{k,\ell}(\lambda)$.
\end{definition}

\begin{lemma}\label{le:G}
For $n$ a non-negative integer and $k$, $\ell$ positive integers we have 
$$
G_{k,\ell}(n)=\sum_{i=1}^{\lfloor\frac{n}{k}\rfloor}p(n-ik\ell) .$$
\end{lemma}
\begin{proof}
Consider the set $\cG_{k,\ell}(n)$ of ordered pairs $(\lambda,i)$ where
$\lambda\vdash n$ is a partition containing a part of
length $ik$ repeated at least $\ell$ times. 
Then $G_{k,\ell}(n)=|\cG_{k,\ell}(n)|$. We have a bijection from
$\cG_{k,\ell}(n)$ to $\displaystyle
\bigcup_{i=1}^{\lfloor\frac{n}{k}\rfloor}\cP(n-ik\ell)$
given by removing $\ell$ parts of length $ik$ from $\lambda$.
\end{proof}

\begin{theorem}\label{th:G}
For $n$ a non-negative integer and $k$, $\ell$ positive integers we have 
$$G_{k,\ell}(n)=F_{k\ell}(n).$$
\end{theorem}
\begin{proof}
This follows from Proposition~\ref{pr:F} and Lemma~\ref{le:G}.
\end{proof}

\begin{remark}
This identity with $k=1$ can be found in Exercise~80 in Chapter~1 of 
Stanley's book~\cite{Stanley:2012a}. 
This is commonly known as Elder's theorem,
and the case $k=1$, $\ell=1$ is known as Stanley's theorem.
In the solution to the exercise in~\cite{Stanley:2012a},
Stanley explains that he discovered the result for $k=1$ and all 
$\ell\ge 1$ in 1972 and 
submitted it to the Problems and Solutions section of 
Amer.\ Math.\  Monthly, where it was rejected with the comment 
``A bit on the easy side, and using only a standard argument.'' 
The statement was reproved by Elder in 1984,
and other proofs can be found in Hoare~\cite{Hoare:1986a} and
Kirdar and Skyrme~\cite{Kirdar/Skyrme:1982a}.
But as pointed out by Gilbert~\cite{Gilbert:2015a}, 
all these proofs were predated by more than a decade
in a 1959 paper of Fine~\cite{Fine:1959a}.
A different generalisation of Elder's theorem was proved by
Andrews and Deutsch~\cite{Andrews/Deutsch:2016a}. 
\end{remark}

\begin{definition} 
For a non-negative integer $n$  and a positive integer $r$ let  
$q(n,r) $ be the  number of partitions of $n$  into $r$ distinct
parts.  
\end{definition}

The following is   well known.

\begin{lemma}\label{l:q(n,r)}  
Let  $r$ be  a positive integer. The generating function for   
$q(n,r) $, $ n\geq 0 $, is given by 
\[  Q_{r}   (t):= \sum_{n=0}^\infty q(n,r)\,t^n = \prod_{i=1}^r\frac{t^i}{1-t^i}. \]
\end{lemma}
\begin{proof}  
By going over to conjugate partitions, one sees that $q(n,r) $ 
is also equal to the number of partitions of  $n$  with  each of 
$1, 2,  \ldots, r$  occurring as  a part and with $r$ as largest
part. Thus, the generating function equals 
\begin{equation*} 
(t + t^2  +  \cdots)  (t^2  + t^4 +  \cdots )    
\cdots   (t^r + t^{2r} +\cdots )  =  \prod_{i=1}^r\frac{t^i}{1-t^i}. 
\qedhere
\end{equation*}
\end{proof}  

\begin{definition}
Let $n$ be a non-negative integer.
If $\lambda\vdash n$ and $k$, $r$ are positive integers,  let 
$c_{k,r}(\lambda)$ be the number of  unordered $r$-tuples of distinct part lengths of $\lambda$
divisible by $k$.
Let  $C_{k, r}(n) = \displaystyle\sum_{\lambda\vdash n} c_{k,r}(\lambda)$  and let
$D_{k, r}(n) = \displaystyle\sum_{\lambda\vdash n}  
\binom{g_{1,k}(\lambda)}{r} $.
\end{definition}

\begin{theorem}\label{t:CandD} 
Let $n$ be a non-negative integer and let $r,k$ be positive integers.
Then
\[C_{k, r}(n) = D_{k,r}(n).\]
The generating function  for   $ C_{k, r}(n) $, $n\geq 0 $  is 
\[  \left(\prod_{i=1}^r \frac{t^{ik}}{1-t^{ik}}\right)P(t).     \]  
\end{theorem} 

\begin{proof}  
Let $\{i_1, \ldots, i_r \}$ be a set of distinct positive numbers. 
The set of partitions of $n$ having $ki_1,\ldots,     ki_r $  
as parts is in bijection  with the set of partitions of  
$n-k( i_1+ \cdots  + i_r)  $  (where as usual  the   set of partitions 
of a negative integer is taken to be the empty set). Similarly, the  
set  of partitions of   $n$   with  each of the parts  
$i_1, \ldots, i_r $ occurring  at least  $k$ times is  in bijection
with the set of partitions of  
$n-k( i_1+ \cdots  + i_r)  $. Thus, the generating functions of both  
$C_{k, r}(n) $ and of  $D_{k,r}(n) $ equals
$Q_r(t^k)P(t)$.
Now the result follows from Lemma~\ref{l:q(n,r)}.
\end{proof} 

\begin{definition}
Let $n$ be a non-negative integer.
If $\lambda\vdash n$, let $e(\lambda) $  denote the  number of  
ordered pairs $(m, m')$  of positive integers  such that $m$ and $m'$ 
are   parts of $\lambda$, $ m \ne m' $,  $\lambda_m \geq 2 $,
and $m' $  is even.
Let  $E(n) = \displaystyle\sum_{\lambda\vdash n} e(\lambda)$.
\end{definition}

\begin{theorem} \label{t:E}  The generating function  for   $ E (n) $, $n\geq 0 $  is 
\[  \frac{t^4+2t^8}{(1-t^2)(1-t^6)}\,P(t).    \]  
\end{theorem} 
\begin{proof}  
Let  $m $  be a positive integer  and let  $ m' =2m''
\ne m $ be an even positive integer. The set of partitions  of  $n$
having  both $m $ and $m'$ as parts    and with  the multiplicity of
$m$ as a part being at least $2$  is in bijection with the
partitions of   $n- 2(m +m'')  $. Thus the  generating function
for  $E(n)  $, $ n \geq 0  $  is  
\[  R (t^2)P(t)   \]   
where   $ R(t)  =\sum_{i \geq 0 } ^{\infty}R_it^i $  and   
$R_i $  equals  the number of  ordered pairs  $ (m, m'')$   of
positive integers  such that   $ i=m+m'' $  and    $m \ne 2 m'' $.
We have 
\[  R_i=\begin{cases}
i-1&  i   \text{ not a multiple of } 3 \\
i-2 &   i \text{ a multiple of } 3 .
\end{cases} \\  \]
Thus, 
\[R(t)  =  \frac{t^2}{(1-t)^2}  -\frac{t^3} {1-t^3} =   \frac{t^2 +2t^4}{ (1-t) (1-t^3)}, \]
proving the result.
\end{proof}

\section{On the Hochschild cohomology of symmetric groups}\label{se:HH}  

If $R$ is a commutative ring of coefficients and $G$ is a finite group, recall 
from~\cite[Theorem~2.11.2]{Benson:1998c} the centraliser 
decomposition of Hochschild cohomology
\[ \HH^*(RG) \cong\! \bigoplus_{\textrm{ccls of }g\in G}\!H^*(C_G(g),R). \]
where the sum is indexed by conjugacy classes of elements of $G$,
and $g$ is a representative of the conjugacy class. For notational
convenience, we adopt the convention that
$S_0$ is the trivial group and that wreath products of the form $G\wr S_0$ are
trivial, for any group $G$. For $m$ a positive integer, we denote by $\bZ/m$ a cyclic
group of order $m$.

\begin{proposition}\label{pr:centraliser}
Let $n$ be a positive integer.
We have an isomorphism of graded vector spaces
\[ \HH^*(\bF_p\sym_n)\cong \bigoplus_{\lambda\vdash n}
\bigotimes_{m=1}^n H^*(\bZ/m\wr\sym_{\lambda_m},\bF_p). \]
\end{proposition}
\begin{proof}
The conjugacy classes of $\sym_n$ correspond to partitions
of $n$. If $\lambda=(n^{\lambda_n}\dots 2^{\lambda_2}1^{\lambda_1})\vdash n$, then an 
element $g$ of cycle type $\lambda$ has as centraliser the direct product of
wreath products
\begin{equation*} 
C_G(g) = \prod_{m=1}^n (\bZ/m \wr \sym_{\lambda_m}). 
\end{equation*}
Applying the K\"unneth formula to this direct product of groups yields the result.
\end{proof}

To compute the dimension of the expression given
in Proposition~\ref{pr:centraliser}, we use Nakaoka's 
description of the cohomology of a wreath product of
finite groups over a field.

\begin{proposition}\label{pr:Nakaoka}
Let $\ell$, $m$ be positive integers. We have an isomorphism of
graded vector spaces
\[ H^*(\bZ/m\wr\sym_{\ell},\bF_p)\cong
H^*(\sym_{\ell},H^*((\bZ/m)^{\ell},\bF_p)). \]
\end{proposition}
\begin{proof}
This follows from Theorem~3.3 of Nakaoka~\cite{Nakaoka:1961a}.
\end{proof}

This is in fact an isomorphism of graded $k$-algebras with the appropriate
algebra structure on the right side, but this will not be needed.

\section{Degree one}\label{se:HH1}  

For the remainder of the paper, we write $M_{\ell}$ for the canonical $\ell$-dimensional
permutation module for $\sym_{\ell}$ over $\bF_p$, where $\ell$ is a positive integer. 
This is induced from the trivial module for $\sym_{\ell-1}$ over $\bF_p$.

\begin{lemma}\label{le:H0}
For $\ell$ a positive integer, we have $H^0(\sym_{\ell},M_{\ell})\cong \bF_p$.
\end{lemma}
\begin{proof}
By Frobenius reciprocity, we have 
$H^0(\sym_{\ell},M_{\ell})\cong H^0(\sym_{\ell-1},\bF_p)$.
\end{proof}

\begin{lemma}\label{le:H1}
Let $\ell$, $m$ be  positive integers.
As modules for $\sym_{\ell}$, we have 
\[ H^1((\bZ/m)^{\ell},\bF_p)
\cong\begin{cases} 0 & p\nmid m,\\ M_{\ell} & p\mid m.
\end{cases} \]
\end{lemma}
\begin{proof}
This follows by looking at the action of $\sym_\ell$ on the two sides of the isomorphism 
\begin{equation*}
H^1((\bZ/m)^{\ell},\bF_p)\cong \Hom((\bZ/m)^{\ell},\bF_p).
\qedhere
\end{equation*}
\end{proof}

\begin{lemma}\label{le:H0,1}
Let $\ell$, $m$ be positive integers.
We have 
\[ \dim_{\bF_p}H^0(\sym_{\ell},H^1((\bZ/m)^{\ell},\bF_p))
=\begin{cases}0&p\nmid m,\\1&p\mid m.\end{cases} \]
\end{lemma}
\begin{proof}
This follows from Lemma~\ref{le:H0} and Lemma~\ref{le:H1}.
\end{proof}

\begin{lemma}\label{le:H0,1sum}
Let $n$ be a positive integer and $\lambda$ a partition of $n$. We have
\[ \sum_{m=1}^n \dim_{\bF_p} H^0(\sym_{\lambda_m}, H^1((\bZ/m)^{\lambda_m},\bF_p))
= g_{p,1}(\lambda) . \]
\end{lemma}

\begin{proof}
Only summands with $\lambda_m\geq 1$ contribute to the left side. 
It follows from Lemma~\ref{le:H0,1} that the left side
counts the number of parts of $\lambda$ 
of distinct lengths divisible by $p$, and this is $g_{p,1}(\lambda)$.
\end{proof}

\begin{lemma}\label{le:H1sum}
Let $n$ be a positive integer and $\lambda$ a partition of $n$. We have 
\[ \sum_{m=1}^n \dim_{\bF_p} H^1(\sym_{\lambda_m}, \bF_p) =
\begin{cases} 
g_{1,2}(\lambda) & p=2,\\
0 & p \ge 3.
\end{cases} \]
\end{lemma}
\begin{proof} 
Only summands with $\lambda_m\ge 1$ contribute to the
left side. For $\ell$ a positive integer
we have $H^1(\sym_\ell, \bF_p)=0$  for $p\ge 3$ 
or $\ell=1$, and  $\dim_{\bF_2} H^1(\sym_\ell, \bF_2) = 1$  if
$\ell\ge 2$.  Thus the left side in the statement is zero for $p\ge 3$.
For $p=2$ this counts the number of distinct part lengths of $\lambda$ 
that appear at least twice, and this is $g_{1,2}(\lambda)$.
\end{proof}

Combining the above results yields the following combinatorial descriptions
of the dimension of $\HH^1(\bF_p\sym_n)$.

\begin{theorem}\label{th:HH1}
Let $n$ be a positive integer.
The dimension of the degree one Hochschild cohomology of the
group algebra of the symmetric group $\sym_n$ on $n$ letters  is given by 
\[ \dim_{\bF_p}\HH^1(\bF_p\sym_n)= 
\begin{cases}
G_{2,1}(n)+G_{1,2}(n) &  p =2, \\
G_{p,1}(n)    &  p \ge 3.
\end{cases}    \]
\end{theorem}

\begin{proof}
Let $\lambda$ be a partition of $n$.
Note that the degree one part of the tensor product on the right side
in Proposition \ref{pr:centraliser} corresponding to the summand indexed
by $\lambda$  is isomorphic to 
\[ \bigoplus_{m=1}^n H^1(\bZ/m\wr\sym_{\lambda_m},\bF_p) . \]
Fix an integer $m$ such that $1\le m\le n$.
Applying Proposition~\ref{pr:Nakaoka} in degree one, we have 
\[ H^1(\bZ/m\wr\sym_{\lambda_m},\bF_p) \cong 
H^0(\sym_{\lambda_m},H^1((\bZ/m)^{\lambda_m},\bF_p)) 
\oplus H^1(\sym_{\lambda_m},\bF_p). \]
Taking the sum over all partitions of $n$ and all $m$ such that $1\le m\le n$
we obtain
\[ \HH^1(\bF_p\sym_n) = \bigoplus_{\lambda\vdash n}\bigoplus_{m=1}^n
\left( H^0(\sym_{\lambda_m},H^1((\bZ/m)^{\lambda_m},\bF_p))
\oplus H^1(\sym_{\lambda_m},\bF_p)\right).  \]
Adding up dimensions, using Lemma \ref{le:H0,1sum} and Lemma
\ref{le:H1sum} yields the result.
\end{proof}

\begin{cor} \label{cor:HH1}
For any positive integer $n$ we have 
\[ \dim_{\bF_p}(\HH^1(\bF_p\sym_n)) =
\begin{cases} 2F_2(n) & p=2, \\
F_p(n) & p \ge 3.
\end{cases} \]
\end{cor}

\begin{proof}
This follows from combining Theorem \ref{th:HH1} and Theorem~\ref{th:G}.
\end{proof}

\begin{proof}[Proof of Theorem~\ref{th:HH}\,\rm (ii)]
Combining Corollary~\ref{cor:HH1}  with Proposition~\ref{pr:F} yields the 
result.
\end{proof}

\section{Degree two}\label{se:HH2} 

\begin{lemma}\label{le:Lambda2M}
Let $\ell$ be a positive integer.
We have 
\begin{align*} 
\dim_{\bF_p}H^0(\sym_{\ell},\Lambda^2(M_{\ell}))&=
\begin{cases}
0 & p\ne 2\text{\rm\ or }\ell = 1,\\
1 & p=2\text{\rm\ and }\ell \ge 2,
\end{cases}\medskip
\\
\dim_{\bF_p}H^1(\sym_{\ell},\Lambda^2(M_{\ell}))&=
\begin{cases}
0 & p\ne 2\text{\rm\ or }\ell = 1,\\
1 & p=2\text{\rm\ and }\ell=2,3,\\
2 & p=2\text{\rm\ and }\ell \ge 4.
\end{cases} 
\end{align*} 
\end{lemma}

\begin{proof}
The module $\Lambda^2(M_{\ell})$ is zero for $\ell =  1$.
For $\ell \ge 2$, this module is
induced from the trivial tensor sign representation $\bF_p\otimes\sgn$ of 
$\sym_{\ell-2}\times\sym_2$, so by Frobenius reciprocity, we have
\[ H^i(\sym_{\ell},\Lambda^2(M_{\ell}))\cong 
H^i(\sym_{\ell-2}\times\sym_2,\bF_p\otimes\sgn). \]
If $p\ne 2$, the module $\bF_p\otimes\sgn$ is not in the principal block, 
so there is no cohomology in any degree. If $p=2$, $\bF_p\otimes\sgn$
is the trivial module, and the computation is straightforward.
\end{proof}

\begin{lemma}\label{le:H2Zm}
Let $\ell$, $m$ be  positive integers.
As an $\bF_p\sym_{\ell}$-module, we have
\[ H^2((\bZ/m)^{\ell},\bF_p) \cong
\begin{cases}
0 & p\nmid m \\
M_{\ell} \oplus \Lambda^2(M_{\ell}) & p \mid m.
\end{cases} \]
\end{lemma}

\begin{proof}
For $p\nmid m$ this is clear, so suppose that $p\mid m$. 
If $p$ is odd, the cohomology is a tensor product of an exterior
algebra on $m$ generators in degree one, tensored with a
polynomial ring on their Bocksteins in degree two.
So in degree two we have the polynomial generators,
which give a copy of $M_{\ell}$, and the products of
two exterior generators, which give a copy of $\Lambda^2(M_{\ell})$.
If $p=2$ and $4\mid m$, the computation is the same. If $p=2$
and $m \equiv 2 \pmod{4}$ then the cohomology is a polynomial 
ring on the degree one generators, so in degree two we have
$S^2(M_{\ell})$. As a module for $\sym_{\ell}$, this
is isomorphic to $M_{\ell} \oplus \Lambda^2(M_{\ell})$.
The squares of degree one generators give a copy of $M_{\ell}$ and the products of
distinct degree one generators give a copy of $\Lambda^2(M_{\ell})$.
\end{proof}

\begin{lemma}\label{le:H0,2}
Let $\ell$, $m$ be  positive integers.
We have
\[ \dim_{\bF_p}
H^0(\sym_{\ell},H^2((\bZ/m)^{\ell},\bF_p))
= \begin{cases}
0 & p\nmid m, \\
1 & 2\ne p \mid m,\\
1 & 2= p \mid m \text{\rm\ and }\ell = 1, \\
2 & 2 = p \mid m \text{\rm\ and } \ell \ge 2.
\end{cases} \]
\end{lemma}

\begin{proof}
By Lemma~\ref{le:H2Zm}, we get zero if $p\nmid m$ 
and two terms when $p\mid m$.
The computation of $H^0(\sym_{\ell},M_{\ell})$ is given
in Lemma~\ref{le:H0}, and the computation of 
$H^0(\sym_{\ell},\Lambda^2(M_{\ell}))$
is given in Lemma~\ref{le:Lambda2M}.
\end{proof}

\begin{lemma}\label{le:H1M}
Let $\ell$ be a positive integer. We have
\[ \dim_{\bF_p}H^1(\sym_{\ell},M_{\ell})=
\begin{cases}
1 & p=2,\ \ell \ge 3,\\
 0 & \text{\rm otherwise.}
\end{cases} \]
\end{lemma}
\begin{proof}
By Frobenius reciprocity, we have
$H^1(\sym_{\ell},M_{\ell})\cong H^1(\sym_{\ell-1},\bF_p)$.
\end{proof}

\begin{lemma}\label{le:H1,1}
Let $\ell$, $m$ be  positive integers. We have 
\[ \dim_{\bF_p}H^1(\sym_{\ell},H^1((\bZ/m)^{\ell},\bF_p)) =
\begin{cases} 
1 & 2=p\mid m \text{\rm\ and } \ell\ge 3,\\
0 & \text{\rm otherwise.}
\end{cases} \]
\end{lemma}

\begin{proof}
By Lemma~\ref{le:H1}, $H^1((\bZ/m)^{\ell},\bF_p)$ is equal to $M_{\ell}$ if $p\mid m$
and zero otherwise. Now use Lemma~\ref{le:H1M}. 
\end{proof}

\begin{lemma}\label{le:H2,0}
Let $\ell$ be a positive integer. We have
\[ \dim_{\bF_2} H^2(\sym_{\ell},\bF_2) =
\begin{cases} 0 & 
\ell=1,\\
1&\ell=2,3,\\ 
2&\ell\ge 4,
\end{cases}  \] 
and $H^2(\sym_{\ell},\bF_p)=0$ for $p\ge 3$.
\end{lemma}

\begin{theorem}\label{th:HH2}
The dimension of the degree two Hochschild cohomology of the 
symmetric groups is given by
\[ \dim_{\bF_p}\HH^2(\bF_p\sym_n) = 
\begin{cases}
G_{2,1}(n)+G_{2,2}(n)+G_{2,3}(n)+G_{1,2}(n)\\
\qquad\qquad\qquad{}+G_{1,4}(n)  +  C_{2,2}(n) +   E(n)+D_{2,2}(n) & p =2,\\
G_{p,1}(n) +   C_{p,2} (n) & p \ge 3.
\end{cases} \]
\end{theorem}

\begin{proof}
By Proposition~\ref{pr:centraliser}, 
\begin{multline} \label{e:kun} \dim_{\bF_p} \HH^2(\bF_p\sym_n)  =  
\sum_{\lambda \vdash n}  \sum_{m=1}^n \dim_{\bF_p}H^2(\bZ/m \wr\sym_{\lambda_m},\bF_p) \\
\qquad{}+   \sum_{\lambda \vdash n} \sum_{\{m, m'\} } 
( \dim _{\bF_p}H^1(\bZ/m\wr\sym_{\lambda_m},\bF_p))  
(\dim _{\bF_p}H^1(\bZ/m'\wr\sym_{\lambda_{m'}},\bF_p))
\end{multline}
where the  indices   $\{m,m'\}$  run over unordered pairs of distinct positive integers  
between $1$ and $n$.

Using Proposition~\ref{pr:Nakaoka}, we have
\begin{multline*} 
H^2(\bZ/m\wr\sym_{\lambda_m},\bF_p) \cong 
H^0(\sym_{\lambda_m},H^2((\bZ/m)^{\lambda_m},\bF_p)) \\
\oplus H^1(\sym_{\lambda_m},H^1((\bZ/m)^{\lambda_m},\bF_p) )\oplus
H^2(\sym_{\lambda_m},\bF_p). 
\end{multline*}
Using Lemmas~\ref{le:H0,2}, \ref{le:H1,1}, and~\ref{le:H2,0}, we get
that  
for $p\ge 3$, the   first term on the right contributes $G_{p,1}(n)$, 
the second and third contribute
zero  to the sum  first sum on the right of Equation~\ref{e:kun}.
For $p=2$, the first term contributes $G_{2,1}(n)+G_{2,2}(n)$, the
second contributes $G_{2,3}(n)$, and
the third contributes $G_{1,2}(n)+G_{1,4}(n)$.

Next we consider the contribution of the  second sum  on the right hand side of
 Equation~\ref{e:kun}.   Suppose  first that  $p\ge 3$.  
As in the proof of  Theorem~\ref{th:HH1}, for any positive integer $m$
 \[  \dim _{\bF_p}H^1(\bZ/m\wr\sym_{\lambda_m},\bF_p)=
 \begin{cases} 
 0 &  p\nmid  m \text{ \rm or } \lambda_m=0  \\
 1 & p\mid  m \text{ \rm  and } \lambda_m \geq 1.  \end{cases}  \] 
Thus  the contribution of  a partition $\lambda\vdash n$  
to the second  sum of Equation~\ref{e:kun} equals $c_{p,2}(\lambda)$ 
and the second sum equals   $C_{p,2} (n) $. This completes the proof
for $p\ge 3$.
Now suppose that $p=2 $. Again as in the  proof of Theorem~\ref{th:HH1} 
 for any positive integer  $m$
  \[\dim _{\bF_p}H^1(\bZ/m\wr\sym_{\lambda_m},\bF_p)=
 \begin{cases} 
 0& \lambda_m =0\\
 0 &  2\nmid  m  \text{ \rm  and }  \lambda_m =1\\
 1 & 2\mid  m  \text{ \rm  and } \lambda_m =1\\
 1 & 2\nmid  m  \text{ \rm  and }  \lambda_m \geq 2\\      
 2   &  2 \mid  m  \text{ \rm  and }  \lambda_m \geq 2.    \end{cases} \]
Thus  the contribution of  a partition $\lambda\vdash n$  to the  second   sum of Equation~\ref{e:kun} equals
\[c_{2,2,}(\lambda)  + e(\lambda)  +  \binom{g_{1,2}(\lambda)} {2}  \]  
and the second sum equals 
\[ C_{2,2}(n) +E(n) +D_{2,2}(n), \]  as required.
\end{proof}

\begin{cor}\label{cor:HH2}
For any positive integer $n$ we have 
\[ \dim_{\bF_p}\HH^2(\bF_p\sym_n)= 
\begin{cases}
2F_2(n)+2F_4(n)+F_6(n)  +  2C_{2,2}(n) +   E(n) & p =2,\\
F_p(n) +   C_{p,2} (n) & p \ge 3.
\end{cases} \]
\end{cor}

\begin{proof}
This follows from combining Theorem~\ref{th:HH2} with Theorem~\ref{th:G}
and the first equality in Theorem~\ref{t:CandD}.
\end{proof}

\begin{proof}[Proof of Theorem~\ref{th:HH}\,\rm (iii)]
Applying Proposition~\ref{pr:F}, Theorem~\ref{t:CandD}   and Theorem~\ref{t:E}
to the terms in the  formula in Corollary \ref{cor:HH2}, 
we obtain
\begin{multline*}
\sum_{n=0}^\infty \dim_{\bF_p}\HH^2(\bF_p\sym_n)\,t^n=\\
\begin{cases}
\displaystyle 
\left(\frac{2t^2}{1-t^2}  +\frac{2t^4}{1-t^4}  + \frac{t^6}{1-t^6}  
+ \frac{2t^6}{(1-t^2) (1-t^4) } + \frac{t^4+2t^8}{(1-t^2) (1-t^6)} \right)P(t)
&p=2,\\[10pt]
\displaystyle \frac{t^p} { (1-t^p)(1-t^{2p})}\,P(t) &p\ge 3.
\end{cases} 
\end{multline*}
Simplifying the expression for $p=2$ completes the proof.
\end{proof}

\section{Proof of Theorem~\ref{th:rational}.} \label{se:HHgeq3} 

\subsection{More combinatorics.} 
If  $f(t) $ and  $g(t)$ are  monic  polynomials  and $\alpha$ is a non-zero constant, then   
we say that  $\alpha $   is  the leading coefficient   of the  rational function      
$  \frac{\alpha f(t) }{g(t) }$  and the difference of the degree of  $f(t)$  and  $g(t)$  is referred to as the total degree of $  \frac{\alpha f(t) }{g(t) }$.
  For a  positive integer  $m$, set   
\[  E_m  (t)   :=   \frac{t^m}{1-  t^m}. \]

We fix a positive  integer  $r$ and let    ${\mathcal V} =\{1,\ldots, r\} $.   We let  
${\mathcal E} $ be  the set of two  element subsets  of  ${\mathcal V}$  and let   
${\mathcal  G}$ be the set  of  undirected simple graphs with vertex set $ {\mathcal V} $.   
Then ${\mathcal G} $    may be identified   with  the  power set of  ${\mathcal E} $  via 
the  map which sends a graph  $Y$ to the subset  of ${\mathcal E} $  consisting of those 
$\{i,j\} $ such that  $i$ and $j$ are connected by an edge in $Y$.  We  refer to the elements 
of ${\mathcal E}$  as edges.

For each $ Y$ in ${\mathcal G} $, let  $e(Y) $   denote  the number  of edges of $Y$   
and let  $c (Y) $  denote the  number of connected components  of $Y$.   Let $C_Y$  
be the set partition of ${\mathcal V}$  induced  by  the equivalence  relation  
corresponding to the  connected components  of $Y$. In other words,   
$i,j\in {\mathcal V} $ are in the same element of $C_Y$  iff  $i$ and $j$  belong to the 
same  connected component of  $Y$.

\begin{lemma}\label{l:graph}  
With the above notation,  $\sum_{Y\in {\mathcal G}  }    (-1)^{c(Y) +e(Y) }    = (-1)^r r!  $.
\end{lemma}  

\begin{proof}    
Let $\Gamma $ be the complete  graph on ${\mathcal V}$.  By  
Tutte~\cite[IX.2.2]{Tutte:1984a}  
 \[\sum_{Y\in {\mathcal G}  }    (-1)^{c(Y) +e(Y) }    =  P_{\Gamma} (-1) \]  where
   $ P_{\Gamma} (x)   $  is   the chromatic polynomial  of $\Gamma $.  
By   \cite[Thm.~IX.26]{Tutte:1984a},   we have
\[  P_{\Gamma} (x) = x (x-1)  \ldots (x - (r-1) )   , \] 
proving the result.
\end{proof}

 Let $k_i,\ell_i $, $1\leq i \leq  r$ be positive integers.  For each non-empty subset    
$A$ of ${\mathcal V} $, set  
 \begin{align*}  k_A &:= \lcm( k_i \colon i \in A)  \\
  \ell_A &:= \sum_{i\in A} \ell_i   \\
 E_A(t) &:=   E_{k_A\ell_A }(t) =  \frac{t^{\ell_Ak_A} }{1-  t^{\ell_Ak_A} }. 
 \end{align*}
 For  a  graph $Y$ in ${\mathcal G} $, set 
\begin{align*}  E_Y(t) &:=    \prod_{A\in C_Y} E_A(t) =
\prod_{A  \in C_Y} \frac{t^{\ell_Ak_A} }{1-  t^{\ell_Ak_A} }. \end{align*}

\begin{lemma}  \label{l:comb} Let $r$,    $k_i $, $\ell_i $, $k_A$, $\ell_A$ be  as above  and let    
$Q(t)  =\sum_{n=1}^{\infty}  q(n)  t^n   $, where 
for  a positive integer $n$,   $q(n)$ equals the number of  ordered  $r$-tuples 
$(u_1,\ldots, u_r) $ of positive integers such that $n=\sum_{i=1}^n \ell_ik_i u_i $ 
and such that $k_iu_i \ne  k_ju_j $ if   $i\ne j $.   Then 
\[   Q(t)    =       \sum_{  Y\in {\mathcal G}  }    (-1)^{e(Y)}       E_Y(t) . \]
Consequently,   $ Q(t)    $ is a rational  function of total degree   $0$ and leading 
coefficient  $(-1)^r r!$.
\end{lemma} 

\begin{proof}  
Let $U(n) $ be the set of ordered tuples  $r$-tuples $(u_1,\ldots, u_r) $ of positive 
integers such that $n=\sum_{i=1}^n \ell_ik_i u_i $.    For  $Y \in {\mathcal G} $, let   
$U_Y (n)   $ be the subset  of   $U (n) $   such that   $ k_iu_i  =   k_ju_j  $  whenever  $i $ 
and $ j$  are in the same connected component of    $Y$.  So, $U(n) = U_{\emptyset}  (n) $  
where $\emptyset $ denotes the  graph  without any edges.  We claim that
\[  q(n)  =  \sum_{Y \in {\mathcal G}  }  (-1)^{e(Y)} |U_Y (n)| .   \]
Indeed,    let $ U'(n)   $   be  the subset   of $U(n)   $ consisting of those $r$-tuples   such 
that  $k_iu_i= k_ju_j $  for some   $i\ne j $.  For  an edge   $e=\{i, j\}   \in  {\mathcal E} $,  
let $U^{e} (n )$  be  the subset   of $U(n)   $ consisting of those $r$-tuples   such that  
$k_iu_i   =    k_ju_j $.    Then  $ U'(n) $ is the union  of the   $ U^{e} (n) $ as  $e$ runs  
over  the  edge  set  $ {\mathcal E}$  and 
  \[ q(n)  =  |U(n)|  -  |U'(n)|.  \]
Now    let $Y  \subseteq   {\mathcal E}$  be a graph.   Then 
$  \bigcap_{e \in Y}  U^{e} (n)   =  U_Y(n)  $. The  claim   now  follows   by
the above equation and   the inclusion-exclusion   principle.

Next,  we  claim that 
\[  \sum_{n=1} ^{\infty}  |U_Y(n) |  t^n =     E_Y(t) .\]   
To see this,   note that  if   $(u_1, \ldots,  u_r ) \in U_Y(n) $, then for  any  $A \in C_Y $  
there exists  a positive number $u'_A  $ such that  $ k_iu_i =  k_A u'_A $ for all $ i\in A$.  
The assignment $ (u_1, \ldots,  u_r )  \mapsto 
 (u'_A)_{A \in    C_Y }   $    induces a  bijection between $U_Y(n)  $  and  the set of    
 tuples $ (u'_A)_{A \in    C_Y } $   of positive  integers  $u'_A$   such that  
 $\sum_{A \in C_Y}   \ell_A k_A   u'_A  =  n$.   The claim follows.

The first assertion is  a consequence of the  two claims.    By
definition,  $E_Y(t) $  is a rational function of total degree $0$
with leading  coefficient    $(-1)^{c(Y)}$.   Taking common denominators
in  $Q(t)$  and applying Lemma~\ref{l:graph}   yields the result. 
\end{proof} 

\noindent
{\bf Remark.}    In  the case $r=1$,  $Q(t)  =
\frac{t^{k_1\ell_1}}{1- t^{k_1\ell_1}} $ is the generating function
of  the sequence $ (G_{k_1, \ell_1} (n) =F_{k_1l_1} (n) )$, $ n\geq 1 $.

\begin{definition}   
Let   $k >0$  and $\ell \geq 0 $  be integers.
Let $\theta_{k, \ell} \colon \bN \times \bN_0  \to
\bZ  $ be defined  by  setting  $\theta_{k, \ell}  (x, y)=1 $   
if $k \mid x $   and  $ y \geq  \ell $ and  $\theta_{k, \ell}  (x, y)  = 0$ otherwise.
\end{definition}

\begin{proposition} \label{p:rational}  
Let $r $  and   $k_i $, $\ell_i $,  $1\leq i \leq  r$  be positive integers  and let   $Q(t)$  
be   as  in  Lemma ~\ref{l:comb}.
Let  
\[G(n)  = \sum_{\lambda \vdash n }\sum_{(m_1, \ldots,  m_r)} 
\prod_{i=1}^r  \theta_{k_i, \ell_i}   (m_i,  \lambda_{m_i}), \]
where the inner sum runs over  ordered  $r$-tuples of    distinct part lengths  
$m_1, \ldots,  m_r $  of $\lambda $.
Then    \[ \sum_{n = 1 }^{\infty}    G(n)  t^n     =   Q(t)  P(t)  .  \]
\end{proposition} 

\begin{proof}   
Let  $u_i $, $ 1\leq i \leq r $   be positive integers  such that  setting 
$m_i := k_iu_i   \ne m_j:=  k_ju_j$,   for any $i\ne j $, $ 1\leq  i, j \leq r $. 
 Then  the  set of  partitions of   
$n $ with  $m_i $ as a part  with multiplicity  at least $\ell_i $, $ i=1,\ldots, r   $  is in 
bijection  with the set of partitions of $n-  (k_1\ell_1u_1 +  \cdots  + k_r\ell_r u_r )$. 
The result follows, using Lemma \ref{l:comb}.
\end{proof}

\subsection{Proof of  Theorem~\ref{th:rational}.}     

\begin{lemma}\label{l:Hwr}  
Let   $d\geq 1  $ and let $h\colon {\mathbb N} \times  {\mathbb N}_0  \to {\mathbb Z} $  
be defined by  
\[ h(m, n) =  \dim_{\bF_p} (H^d(\bZ/m\wr\sym_{n},\bF_p)),   \    m \in
  {\mathbb N},   \  n \in   {\mathbb N}_0.  \] 
Then  $h $ is  a   ${\mathbb Z}$-linear  combination  of   $\theta_{k_i,  \ell_i} $   for some 
 positive integers  $k_i , \ell_i  $,  $ 1\leq i \leq s $, with the property that if $\theta_{k_i,\ell_i}$
 appears with a non-zero coefficient, then $k_i\in \{1,p\}$ when $p$ is odd,
 and $k_i\in \{1,2,4\}$ when $p=2$, for all $i$ such that $1\leq i\leq s$.
\end{lemma}  

\begin{proof}   
For  a positive integer  $k$, let   $\epsilon_k\colon
{\mathbb N}  \to {\mathbb Z} $   be  the function defined by
$\epsilon_k(m)  =1 $ if $ k \mid   m$ and $\epsilon_k(m)  =0 $
otherwise.   For a  non-negative integer   $\ell $, let $\gamma_\ell
\colon {\mathbb N}_0 \to {\mathbb Z} $   be defined by 
$\gamma_\ell(n) =1 $ if $  n \geq  \ell $ and $\gamma_{\ell} (n)    =0 $  otherwise.
So,   $\theta_{k, \ell} (m,n)  = \epsilon_k( m) \gamma_{\ell} (n) $
for all $ k, m >0$ and  all $\ell, n \geq 0 $.

Suppose that $p$ is odd. By  Proposition~\ref{pr:Nakaoka} and the
structure of the  cohomology    of cyclic groups,  $h(m, -) =  h(p, -)
$ if   $p \mid m $   and    $h(m, -) =  h(1, -)   $  if $p \nmid m $. Thus,
for all $ n\geq 0 $,  
\[ h(m,n)   = \epsilon_1 (m)   h(1, n)  -\epsilon_p (m)   h(1, n)  +  \epsilon_p(m) h(p,n).    \]  
Similarly, if $p=2 $, then $h(m, -) =  h(4, -)   $ if   $4\mid m $,
$h(m, -) =  h(2, -)   $  if  $m\equiv 2 \pmod 4$ 
and  $h(m, -) =  h(1, -)   $  if     $ m$  is odd.   Hence,  
\[  h(m,n) =\epsilon_1(m) h(1, n)  -\epsilon_2 (m)  h(1, n)  +
  \epsilon_2 (m)   h(2,n)  - \epsilon_4 (m) h(2,n)   + \epsilon_4 (m)
  h(4,n). \] 

Let  $k \in \{1,p ,4\}  $.    By  \cite[Prop.~1.6]{Hatcher/Wahl:2010a}, the sequence
$h (k,n ),  n \geq 0 $    eventually stabilises, that  is, there
exists  $s_k >0 $    such that  $h (k,n )= h(k, s_k)  $ for all $ n
\geq  s_k $. 
Here, we note that  for $k=1 $,  this is Nakaoka's result on  the
stability of  the  cohomology of   symmetric groups.  
Since   $d \geq 1 $,  $h(k, 0)= 0 $.   It follows  that    there exist
integers   $u_{i,k}$, $ 1  \leq i  \leq n_k $  such that  
\[  h(k, n)  = \sum_{\ell =1}^ {s_k}  u_{i,k}   \gamma_{\ell }  (n)    \]   
for all $n \geq 0 $.
Combining this  with the  previous displayed equations  yields the
result. 
\end{proof}  

Theorem~\ref{th:rational}    is a  consequence of  the above lemma  and 
Lemma~\ref{l:comb}  as we now  show.   For each $d \geq 1 $, we fix a   
${\mathbb Z}$-linear combination of  $\theta_{k,l} $  representing the
function $h$   as in Lemma~\ref{l:Hwr}.

\begin{proof}[{Proof of Theorem \ref{th:rational}}]     
For positive integers  $d$ and $n$  let $A_{d,n} $ be the set of  ordered tuples 
$(d_1, \ldots, d_n)  $  of   non-negative integers  $d_m$ such that 
$\sum_{m=1} ^n d_m = d $, considered as $\sym_n $-set via place permutation. Let  
$A'(d,n)$   be the   ($\sym_n$-invariant) subset of $A_{d,n}  $ consisting of tuples  in 
which  each component is strictly positive.

 By Proposition~\ref{pr:centraliser},
\[\dim_{\bF_p} ( \HH^d(\bF_p\sym_n)  )  =   \sum_{\lambda\vdash n}
\sum_{(d_1, \ldots , d_n)  \in A_{d,n} }    \prod_{m=1}^n \dim_{\bF_p}
( H^{d_m} (\bZ/m\wr\sym_{\lambda_{m} },\bF_p) ). \]
Thus  in order to prove   Theorem~\ref{th:rational} it suffices to show  that for any   
$\sym_n $-orbit ${\mathcal O}$ 
of $A_{d,n} $ the  generating function of the sequence $(G_n)_{n \ge 1}    $   where  
\[  G_n  :=   \sum_{\lambda\vdash n}
\sum_{(d_1, \ldots , d_n)  \in  {\mathcal O}}    \prod_{m=1}^n 
\dim_{\bF_p} ( H^{d_m} (\bZ/m\wr\sym_{\lambda_{m} },\bF_p) ) \]
is of the form $ S(t)  P(t) $  for some rational  function  $S(t) $.

Let ${\mathcal O} $ be as above and let  ${\mathcal O}_0$ be the  subset  of 
${\mathcal O}$ consisting  of the tuples  $(d_1, \ldots, d_n)$  such  that  if   
$d_i \geq 1 $   then   $\lambda_i > 0 $. Since $H^i( Z/m \wr S_0,  \bF_p) = 0 $ for  
all $ m\ge 1$ and all $i \ge 1$,   only elements of ${\mathcal O}_0 $ contribute to $G_n $.     
Now   for  each  $r >0$,     and each   size $r$  subset $M$ of   the set  of distinct part 
lengths  of  $\lambda$,  let     $ {\mathcal O}_{0,M} $   be the subset of   ${\mathcal O}_0$ 
consisting of those tuples $(d_1, \ldots, d_m) $   such that  $d_m >0 $ if and only 
if $m \in M$.  Then  ${\mathcal O} _{0,M} $   is in bijection with  a union of   
$\sym_r$-orbits    of $A'_{d, r}  $  via the  map which sends   an $n$-tuple 
$ (d_1, \ldots, d_n) $ to the  subtuple  of   positive $d_m $'s.     Further, note that  
$H^0(G, \bF_p)  $ is one-dimensional for  all finite groups $G$. Hence, for any 
$(d_1, \ldots, d_n)  \in {\mathcal O}$, the   product  $ \prod_{m=1}^n 
\dim_{\bF_p} ( H^{d_m} (\bZ/m\wr\sym_{\lambda_{m} },\bF_p) )  $ equals  the 
subproduct  ranging over the   $m$'s such that $d_m> 0 $.    Thus, running over  all   
size $r$ subsets of  part lengths of  $\lambda$  and over  all $r >0$,  it follows that    
$G_n $  may be  written as
\[ \sum_{\lambda\vdash n}
\sum_{(e_1, \ldots , e_r)  \in  {\mathcal O'}}    \sum_{(m_1 <  \cdots <  m_r)}  
\prod_{i=1} ^r \dim_{\bF_p} ( H^{e_i} (\bZ/m_i\wr\sym_{\lambda_{m_i} },\bF_p) ) \]
where  $ {\mathcal O'}$ is an   $\sym_r$-orbit  of $A'_{d, r}  $ for some $r$, and where 
in the inner sum   $(m_1  <  \cdots <  m_r) $ runs over   collections  of  $r$ distinct    
part lengths   of $\lambda $.

By  the above  and   Lemma~\ref{l:Hwr}     $G_n $  is a sum  of   terms of the form
\[   G'_n= \sum_{\lambda\vdash n}
 \sum_{  ((a_1, b_1), \ldots,  (a_r, b_r))  \in  {\mathcal X}}   \sum_{(m_1 <  \cdots <  m_r)}   
 \theta_{a_1, b_1} (m_1, \lambda_1)   \ldots  \theta_{a_r, b_r}      (m_r, \lambda_r) \]
where    ${\mathcal X}  $ is an $\sym_r$-orbit of  $r$-tuples of the  form   
$((a_1, b_1), \ldots  (a_r, b_r))$ with $\sym_r $ acting again by place permutations. 
It therefore  suffices to  show that the generating function  of $(G'_n) $, $ n\geq 1 $ is   
the product of  a rational function  with  $P(t)$.  But this generating  function  is   
$  \frac{1}{c}     Q(t)$, where $Q(t) $ is as in Lemma~\ref{l:comb}  and where   
$c  =   \frac{r!}  { |\mathcal X|}  $  is the order of the  $\sym_r$-stabiliser 
of  an element $ ((a_1, b_1), \ldots,  (a_r, b_r)) $   of ${\mathcal  X}$. 
\end{proof}

\bibliographystyle{amsplain}
\bibliography{../repcoh}

\end{document}